\documentclass[12pt,a4paper]{article}

\usepackage{amsmath, amssymb, amsthm}
\usepackage[dvips,final]{epsfig}
\usepackage[dvips,final]{graphicx}

\usepackage{color}

\definecolor{green1}{rgb}{0.1,0.65,0}
\definecolor{viol}{rgb}{0.8,0,0.2}
\definecolor{cyan1}{rgb}{0,0.7,0.3}

\def\e{\varepsilon}
\def\dd{\,\mathrm{d}}
\def\dive{\mathrm{\,div\,}}
\def\diam{\mathrm{diam}}

\def\Lin{\mathrm{Lin}}

\def\bfj{\mathbf{1}}
\def\real{\mathbb{R}}
\def\nat{\mathbb{N}}
\def\tens{\real^{3\times 3}_{\mathrm{sym}}}

\def\AA{\mathbf{A}}
\def\BB{\mathbf{B}}
\def\Bcal{\mathcal{B}}
\def\Ecal{\mathcal{E}}
\def\Wcal{\mathcal{W}}
\def\play{\mathfrak{p}}

\def\nas{\nabla_s}

\def\io{\int_{\Omega}}
\def\ipo{\int_{\partial\Omega}}
\def\oi{^{(i)}}
\def\om{^{(n)}}

\newcommand{\for}{\ \mbox{for } \,}
\newcommand{\ale}{\ \mbox{a.\,e. }}

\def\sip{\sigma^p}
\def\scal#1{\left\langle #1\right\rangle_{\AA^p}}

\def\proz{\mathrm{Proj}_Z}

\def\be{\begin{equation}\label}
\def\ee{\end{equation}}

\def\ber{\begin{eqnarray}}
\def\eer{\end{eqnarray}}

\def\bers{\begin{eqnarray*}}
\def\eers{\end{eqnarray*}}

\newfont{\ctv}{msam10}

\newcommand{\bbox}{\mbox{\ctv \symbol{4}}}
\def\QED{{${}\hfill\bbox$}}
\newenvironment{pf}[1]{\par\vskip1mm{\noindent\it #1.}\ }{\QED\par
\vskip2mm}

\def\bpf{\begin{pf}}
\def\epf{\end{pf}}

\newtheorem{theorem}{Theorem}[section]
\newtheorem{lemma}[theorem]{Lemma}

\newtheorem{hypothesis}[theorem]{Hypothesis}
\newtheorem{proposition}[theorem]{Proposition}

\numberwithin{equation}{section}

\begin{document}

\title{Solvability of an unsaturated porous media flow problem with thermomechanical interaction
\thanks{Supported by the GA\v CR Grant GA15-12227S and RVO: 67985840 (PK),
and by the FP7-IDEAS-ERC-StG Grant \#256872 EntroPhase 
and GNAMPA (Gruppo Nazionale per l'Analisi Matematica, la Probabilit\`a
e le loro Applicazioni) of INdAM (Istituto Nazionale di Alta
Matematica) (ER).}}

\author{Bettina Albers
\thanks{University of Duisburg-Essen, 
Geotechnical Engineering,
D-45117 Essen, Germany, E-mail {\tt bettina.albers@uni-due.de}.}
\and Pavel Krej\v c\'{\i}
\thanks{Institute of Mathematics, Czech Academy of Sciences, \v Zitn\'a~25,
CZ-11567~Praha 1, Czech Republic, E-mail {\tt krejci@math.cas.cz}.}
\and Elisabetta Rocca
\thanks{WIAS Berlin, Mohrenstr.~39, D-10117 Berlin, Germany, E-mail {\tt rocca@wias-berlin.de}.}
}

\maketitle


\begin{abstract}
A PDE system consisting of the momentum balance, mass balance, and energy
balance equations for displacement, capillary pressure, and temperature as a model for
unsaturated fluid flow in a porous viscoelastoplastic solid is shown to admit
a solution under appropriate assumptions on the constitutive behavior.
The problem involves two hysteresis operators accounting for plastic and capillary hysteresis.
\end{abstract}


\section*{Introduction}\label{int}

In a deformable porous solid filled with two immiscible fluids (water and air, say),
two sources of hysteresis
are observed: the solid itself is subject to irreversible plastic deformations, and
the fluid flow exhibits capillary hysteresis which is often explained by the surface tension
on the interfaces between the two fluids. A lot of works have been devoted to this phenomenon,
see, e.\,g., \cite{AM,Murphys,CEER,Hab,Flynn,flpokr,hr}. Mathematical analysis of various
mechanical porous media models with capillary hysteresis and without temperature effects
has been carried out in \cite{bv1,bv2,show,shst}. A PDE system
for elastoplastic porous media flow with thermal interaction was derived in \cite{ak},
but the existence of solutions was only proved for the isothermal case.

Here, we focus on the qualitative analysis of the model derived in \cite{ak}, assuming in
addition that the heat conductivity depends in a controlled way on the temperature
similarly as in the phase transition model considered in \cite{rr}.
Indeed, we borrow here some techniques employed in 
\cite{rr} and \cite{ak}  in order to prove existence of
a weak solution  for the initial boundary value problem 
associated with the PDE system coupling the momentum balance
(cf.~\eqref{e23}), featuring, in particular, a thermal expansion term depending on the 
temperature field, with a mass balance  (cf.~\eqref{e24})
ruling the evolution of the capillary pressure, and an energy balance (cf.~\eqref{e25})
displaying, in particular, quadratic dissipative terms on the right hand side.

The main mathematical difficulties are related
to the low regularity of the temperature field mainly due to the presence of the 
high order dissipative terms in the internal energy balance.
This is the reason why we need to employ here a key-estimate (cf.~\eqref{est2}), 
already exploited in \cite{fpr} and 
more recently in \cite{rr} for the analysis of non-isothermal phase transition models.
Roughly speaking, since the test of the internal energy balance by the 
temperature $\theta$ is not allowed, we test by a suitable negative power of it and use the growth condition of the heat conductivity $\kappa$ in Hypothesis \ref{h1}~(ii). Another key point in our proof is the $L^\infty$ estimate we get on the pressure which entails a bound in a proper negative Sobolev space for the time derivative of the absolute temperature, which turns out to be another fundamental ingredient in order to pass to the limit in our approximation scheme.

The structure of the paper is as follows. The model from \cite{ak} is briefly summarized
in Section \ref{mod}. In Section \ref{hys} we recall the definitions
and main results of the theory of hysteresis operators that are used here. Section \ref{mai}
contains the mathematical hypotheses and statements of the main results. In Section \ref{reg}
we regularize the problem by adding a small parameter $\delta$ accounting for ``micro-movements''
and a large cut-off parameter $R$ to control the nonlinearities, and solve the regularized problem by
the standard Faedo-Galerkin method. In Section \ref{pt1} we let $\delta$ tend to $0$ and $R$ to $\infty$
and prove that in the limit, we obtain a solution to the original problem.


\section{The model}\label{mod}

Consider a domain $\Omega\subset \real^3$ filled with a deformable solid matrix material with pores
containing a mixture of liquid and gas. We state the balance laws in referential (Lagrangian)
coordinates, assume the deformations small, and denote for $x\in \Omega$ and time $t\in [0,T]$
\begin{description}
\item $u(x,t)$ ... displacement vector of the referential particle $x$ at time $t$;
\item $\e(x,t) = \nas u(x,t)$ ... linear strain tensor, $(\nas u)_{ij} :=
\frac12 \left(\frac{\partial u_i}{\partial x_j} + \frac{\partial u_j}{\partial x_i}\right)$;
\item $\sigma(x,t)$ ... stress tensor;
\item $p(x,t)$ ... capillary pressure;
\item $\theta(x,t)$ ... absolute temperature;
\item $A(x,t)$ ... relative gas content.
\end{description}

\begin{figure}[htb]
\centerline{\psfig{file=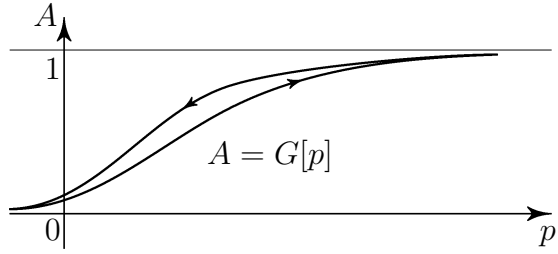,width=7.5cm}}
\caption{Pressure-saturation hysteresis diagram}\label{f3}
\end{figure}

For the stress $\sigma$ and gas content $A$ we assume the empirical constitutive relations
\ber\label{e2}
\sigma &=&  \BB\e_t + P[\e] + (p - \beta(\theta - \theta_c))\bfj\,,\\ \label{e3}
A &=& G[p]\,,
\eer
where $P$ is a hysteresis operator describing the elastoplastic response of the solid,
see Subsection \ref{plas},
$\BB$ is a constant symmetric positive definite fourth order viscosity tensor, $\beta\in \real$ is the
relative solid-liquid thermal expansion coefficient, $\theta_c>0$ is a fixed referential temperature,
$\bfj$ is the Kronecker tensor,
and $G$ is a hysteresis operator as on Figure \ref{f3}, see Subsection \ref{prei}.
We will see that both hysteresis operators $P$ and $G$ admit hysteresis potentials $V_P$ (clockwise)
and $V_G$ (counterclockwise) and dissipation operators $D_P, D_G$ such that for all absolutely
continuous inputs $\e, p$, the inequalities
\be{e10}
P[\e]:\e_t - V_P[\e]_t = \|D_P[\e]_t\|_*\,, \quad G[p]_t p - V_G[p]_t = |D_G[p]_t|
\ee
hold almost everywhere, where $\|\cdot\|_*$ is a seminorm in $\tens$.

We assume the heat conductivity $\kappa(\theta)$ depending on $\theta$, and as in  
\cite{ak}, we obtain the system of momentum balance \eqref{e23}, mass balance \eqref{e24}, and
energy balance equations \eqref{e25} in the form
\ber\label{e23}
\rho_S u_{tt} \!\!&{\!=\!}&\!\!\dive (\BB\nas u_t + P[\nas u]) + \nabla p - \beta \nabla\theta + g\,,\\ \label{e24}
G[p]_t \!\!&{\!=\!}&\!\! \dive u_t + \frac{1}{\rho_L} \dive(\mu(p)\nabla p)\,,\\ \nonumber
c_0 \theta_t \!\!&{\!=\!}&\!\! \dive(\kappa(\theta)\nabla\theta)
+ \|D_P[\nas u]_t\|_* + |D_G[p]_t| +  \BB\nas u_t:\nas u_t
+ \frac{1}{\rho_L} \mu(p)|\nabla p|^2\\ \label{e25}
&&-\beta \theta \dive u_t\,,\qquad
\eer
where $c_0>0$ is a constant specific heat, $\rho_S, \rho_L$ are the mass densities of the solid and liquid, respectively, $\BB$
is a positive definite viscosity matrix, $\beta \in \real$ is the relative thermal expansion coefficient,
and $g$ is a given volume force (gravity, e.g.).

We complement the system with initial conditions
\be{ini}
u(x,0) = u^0(x)\,,\quad u_t(x,0) = u^1(x)\,,\quad p(x,0) = p^0(x)\,, \quad
\theta(x,0) = \theta^0(x)\,, 
\ee
and boundary conditions
\be{bou}
\left.
\begin{array}{rcl}
u &=& 0\\
\mu(p)\nabla p\cdot n &=& \gamma_p(x)(p^* - p)\\
\kappa(\theta) \nabla\theta\cdot n &=& \gamma_\theta(x)(\theta^* - \theta)
\end{array}
\right\} \mbox{on }\ \partial\Omega\,,
\ee
where $\gamma_p, \gamma_\theta: \partial\Omega \to [0, \infty)$ are given smooth functions.



\section{Hysteresis operators}\label{hys}

We recall here the basic concepts of the theory of hysteresis operators that are needed in the sequel.

\subsection{The operator $P$}\label{plas}

In \eqref{e2}, $P$ stands for the elastoplastic part $\sigma^{ep}$ of the stress tensor $\sigma$.
We proceed as in \cite{lch} and assume that $\sigma^{ep}$ can be represented as the sum
$\sigma^{ep} = \sigma^e + \sip$ of an elastic component $\sigma^e$ and plastic component $\sip$.
While $\sigma^e$ obeys the classical linear elasticity law
\be{he1}
\sigma^e = \AA^e\e
\ee
with a constant symmetric positive definite fourth order elasticity tensor $\AA^e$,
for the description of the behavior of $\sip$, we split also the strain tensor $\e$ into the sum
$\e = \e^e + \e^p$ of the elastic strain $\e^e$ and plastic strain $\e^p$, and assume
\be{he2}
\sip = \AA^p\e^e
\ee
again with a constant symmetric positive definite fourth order elasticity tensor $\AA^p$,
and for a given time evolution $\e(t)$ of the strain tensor, $t \in [0,T]$, we require $\sip$
to satisfy the constraint
\be{h3}
\sip(t) \in Z \quad \forall t\in [0,T]\,,
\ee
where $Z\subset \tens$ is the domain of admissible plastic stress components.
We assume that it has the form
\be{z0}
Z = Z_0 \oplus \Lin\{\bfj\}\,,
\ee
where $\Lin\{\bfj\}$ is the 1D space spanned by the Kronecker tensor $\bfj$
and $Z_0$ is a bounded convex closed subset with $0$ in its interior
of the orthogonal complement $\Lin\{\bfj\}^\perp$
of $\Lin\{\bfj\}$ (the {\em deviatoric space\/}). 
The boundary $\partial Z$ of $Z$ is the {\em yield surface\/}. The time evolution
of $\e^p$ is governed by the {\em flow rule\/}
\be{h4}
\e^p_t : (\sip - \tilde\sigma) \ge 0 \quad \forall \tilde\sigma \in Z\,,
\ee
which implies that
\be{h4a}
\e^p_t : \sip =  M_{Z^*}(\e^p_t)\,,
\ee
where $M_{Z^*}$ is the Minkowski functional of the polar set $Z^*$ to $Z$.
The physical interpretation of \eqref{h4} is the {\em maximal dissipation principle\/}.
Geometrically, it states that the plastic strain rate $\e^p_t$ points in the outward normal
direction to the yield surface at the point $\sip$. Indeed, if $\sip$ is in the interior
of $Z$, then $\e^p_t = 0$.

It follows from \eqref{z0} that there exist $\sip_0 \in Z_0$ (plastic stress deviator)
and $c\in \real$ (pressure) such that
$\sip = \sip_0 - c\bfj$. On the other hand, putting in \eqref{h4} $\tilde\sigma = \sip_0 - \rho\bfj$
for an arbitrary $\rho \in \real$, we obtain $\e^p_t : \bfj = 0$ (in other words,
no volume changes occur during plastic deformation), so that
\be{h4b}
M_{Z^*}(\e^p_t) = \e^p_t : \sip_0 \le \diam(Z_0)|\e^p_t|.
\ee
We can eliminate the internal variables $\e^e, \e^p$ and write \eqref{h4} in the form
\be{h6}
(\e_t - (\AA^p)^{-1}\sip_t) : (\sip - \tilde\sigma) \ge 0 \quad \forall \tilde\sigma \in Z\,.
\ee
We now define a new scalar product $\scal{\cdot,\cdot}$ in $\tens$ by the formula $\scal{\xi,\eta}
= (\AA^p)^{-1}\xi:\eta$ for generic tensors $\xi, \eta$, and rewrite \eqref{h6} as
\be{h6a}
\scal{\AA^p\e_t - \sip_t, \sip - \tilde\sigma} \ge 0 \quad \forall \tilde\sigma \in Z\,.
\ee
We prescribe a canonical initial condition for $\sip$, namely
\be{h5}
\sip(0) = \proz(\AA^p\e(0))\,,
\ee
where $\proz$ is the orthogonal projection $\tens \to Z$ with respect to the scalar product
$\scal{\cdot,\cdot}$ and is characterized by the variational inequality
\be{h7}
x = \proz(u) \ \Longleftrightarrow \ x\in Z\,, \ \scal{u-x, x-y} \ge 0 \quad \forall y \in Z\,.
\ee
We list here some properties of the variational problem \eqref{h3}, \eqref{h6a}, \eqref{h5}.
The proof can be found in \cite[Chapter I]{book}.

\begin{proposition}\label{p3}
For every $\e \in W^{1,1}(0,T;\tens)$ there exists a unique $\sip \in W^{1,1}(0,T;\tens)$
satisfying \eqref{h3}, \eqref{h6a}, \eqref{h5}. The solution mapping
$$
P_0 : W^{1,1}(0,T;\tens)\to W^{1,1}(0,T;\tens): \e \mapsto \sip
$$
has the following properties.
\begin{itemize}
\item[{\rm (i)}] For all $\e\in W^{1,1}(0,T;\tens)$ we have $|P_0[\e]_t| \le |\e_t|$ a.\,e.,
$P_0 : W^{1,1}(0,T;\tens)\to W^{1,1}(0,T;\tens)$ is strongly continuous, and admits
an extension to a strongly continuous mapping $C([0,T];\tens)\to C([0,T];\tens)$.
\item[{\rm (ii)}] There exists a constant $C>0$ such that for every $\e_1, \e_2 \in W^{1,1}(0,T;\tens)$
and every $t \in [0,T]$ we have
\be{h8}
|P_0[\e_1](t) - P_0[\e_2](t)| \le C\left(|\e_1(0) - \e_2(0)|
+ \int_0^t|(\e_1)_t(\tau) - (\e_2)_t(\tau)|\dd\tau\right).
\ee
\item[{\rm (iii)}] For all $\e\in W^{1,1}(0,T;\tens)$, the energy balance equation
\be{enerp}
P_0[\e]:\e_t - \frac12 \left((\AA^p)^{-1}P_0[\e]:P_0[\e]\right)_t = M_{Z^*}\left((\e - (\AA^p)^{-1}P_0[\e])_t\right)
\ee
is satisfied almost everywhere in $(0,T)$, where $M_{Z^*}$ is the Minkowski functional of the polar set
$Z^*$ to $Z$.
\end{itemize}
\end{proposition}

It follows from Proposition \ref{p3} that the operator $P$ in \eqref{e2} can be represented in the form
\be{P}
P[\e] = \AA^e \e + P_0[\e]\,,
\ee
and the first energy identity in \eqref{e10} holds with the choice
\be{VP}
V_P[\e] = \frac12 \AA^e \e:\e + \frac12 (\AA^p)^{-1}P_0[\e]:P_0[\e]\,, \quad 
D_P[\e] = \e - (\AA^p)^{-1}P_0[\e]\,, \quad \|\cdot\|_* = M_{Z^*}(\cdot)\,.
\ee


\subsection{The operator $G$}\label{prei}

Similarly as in \eqref{P}, the operator $G$ is considered as a sum
\be{G}
G[p] = f(p) + G_0[p]\,,
\ee
where $f$ is a monotone function satisfying Hypothesis \ref{h1}\,(iii) below, and $G_0$
is a Preisach operator that we briefly describe here.

The construction of the Preisach operator $G_0$ is also based on a variational inequality of the type
\eqref{h6a}. More precisely, for a given input function $p \in W^{1,1}(0,T)$ and a memory parameter
$r>0$, we define the function $\xi_r(t)$ as the solution of the variational inequality
\be{play}
\left\{
\begin{array}{ll}
|p(t) - \xi_r(t)| \le r & \forall t\in [0,T]\,,\\
(\xi_r)_t(p(t) - \xi_r(t) - z) \ge 0 & \mbox{a.~e.}\ \forall z \in [-r,r]\,,
\end{array}
\right.
\ee
with a prescribed initial condition $\xi_r(0) \in [p(0) - r, p(0) + r]$. 

This is indeed a scalar version of \eqref{h6a} with $Z$ replaced by the interval $[-r,r]$,
$\e$ replaced by $p$, and $\sip$ replaced by $p - \xi_r$. Here, we consider the whole continuous family
of variational inequalities \eqref{play} parameterized by $r>0$. We introduce the memory state space
\be{lam}
\Lambda = \{\lambda \in W^{1,\infty}(0,\infty): |\lambda'(r)| \le 1 \mbox{ a.~e.}\}\,,
\ee
and its subspace
\be{lamk}
\Lambda_K = \{\lambda \in \Lambda: \lambda(r) = 0 \mbox{ for } r \ge K\}\,.
\ee
We fix $K>0$ and an initial state $\lambda_{-1} \in \lambda_K$, and choose the initial condition as
\be{iniplay}
\xi_r(0) = \max\{p(0) - r, \min\{\lambda_{-1}(r), p(0) + r\}\}\,.
\ee
We have indeed for all $r>0$ the natural bound
\be{binipl}
|\xi_r(0)| \le \max\{|p(0)|, K\}\,.
\ee
The mapping $\play_r: W^{1,1}(0,T) \to W^{1,1}(0,T)$ which with each $p\in W^{1,1}(0,T)$ associates
the solution $\xi_r = \play_r[p] \in W^{1,1}(0,T)$ of \eqref{play}, \eqref{iniplay} is called
the {\em play\/}. This concept goes back to \cite{kp}, and the proof of the following statements
can be found, e.~g., in \cite[Chapter II]{book}.

\begin{proposition}\label{p5}
For each $r>0$, the mapping $\play_r: W^{1,1}(0,T) \to W^{1,1}(0,T)$ is Lipschitz continuous
and admits a Lipschitz continuous extension to $\play_r: C[0,T] \to C[0,T]$ in the sense that
for every $p_1, p_2 \in C[0,T]$ and every $t \in [0,T]$ we have
\be{lipc}
|\play_r[p_1](t) - \play_r[p_2](t)| \le \max_{\tau \in [0,t]}|p_1(\tau) - p_2(\tau)|\,.
\ee
Moreover, for each $p \in W^{1,1}(0,T)$, the energy balance equation
\be{enerpl}
\play_r[p]_t p - \frac12\left(\play_r^2[p]\right)_t = \left|r \play_r[p]_t\right|
\ee
and the identity
\be{mono}
\play_r[p]_t p_t = (\play_r[p]_t)^2
\ee
hold almost everywhere in $(0,T)$.
\end{proposition}

\begin{proposition}\label{p6}
Let $\lambda_{-1} \in \Lambda_K$ be given, and let $\{\play_r: r>0\}$ be the family of play operators.
Then for every $p \in C[0,T]$ and every $t \in [0,T]$ we have
\begin{itemize}
\item[{\rm (i)}] $\play_r[p](t) = 0$ for $r \ge K^*(t):= \max\{K, \max_{\tau \in [0,t]}|p(\tau)|\}$;
\item[{\rm (ii)}] The function $r \mapsto \play_r[p](t)$ belongs to $\Lambda_{K^*(t)}$.
\end{itemize}
\end{proposition}

Given a nonnegative function $\rho \in L^1((0,\infty)\times \real)$ (the {\em Preisach density\/}),
we define the Preisach operator $G_0$ as a mapping that with each $p \in C[0,T]$ associates the integral
\be{pre}
G_0[p](t) = \int_0^\infty \int_0^{\play_r[p](t)} \rho (r,v)\dd v\dd r\,.
\ee
For our purposes, we adopt the following hypothesis on the Preisach density.

\begin{hypothesis}\label{h2}
There exists a function $\rho^* \in L^1(0,\infty)$ such that for a.~e. $v \in \real$ we have
$0 \le \rho(r,v) \le \rho^*(r)$, and we put
\be{crho}
C_\rho = \int_0^\infty\int_{-\infty}^\infty\rho(r,v)\dd v\dd r\,, \quad
C_\rho^* = \int_0^\infty\rho^*(r)\dd r\,.
\ee
\end{hypothesis}

For the reader who is more familiar with the original Preisach construction in \cite{pr}
based on non-ideal relays, let us just point out that for integrable densities, the variational setting
in \eqref{pre} is equivalent, as shown in \cite{max}.

{}From \eqref{enerpl}, \eqref{mono}, and \eqref{pre} we immediately deduce the Preisach
energy identity
\be{enerpr}
G_0[p]_t p - V_0[p]_t = |D_0[p]_t| \ \mbox{ a.~e.}
\ee
provided we define the Preisach potential $V_0$ and the dissipation operator $D_0$ by the integrals
\be{prpot}
V_0[p](t) = \int_0^\infty \int_0^{\play_r[p](t)} v \rho (r,v)\dd v\dd r\,, \quad
D_0[p](t) = \int_0^\infty \int_0^{\play_r[p](t)} r \rho (r,v)\dd v\dd r\,.
\ee
The second identity in \eqref{e10} then holds with the choice
\be{VG}
V_G[p] = p f(p) - \int_0^p f(z)\dd z + V_0[p]\,, \quad D_G[p] = D_0[p]\,.
\ee
A straightforward computation shows that $G_0$ (and, consequently, $G$) are Lipschitz continuous
in $C[0,T]$. Indeed, using \eqref{lipc} and Hypothesis \ref{h2}, we obtain for $p_1, p_2 \in C[0,T]$ and
$t \in [0,T]$ that
\be{lipg}
|G_0[p_2](t) - G_0[p_1](t)| = \left|\int_0^\infty\int_{\play_r[p_1](t)}^{\play_r[p_2](t)}\rho(v,r)\dd v\dd r\right|
\le C_\rho^* \max_{\tau \in [0,t]} |p_2(\tau) - p_1(\tau)|\,.
\ee
We similarly get, using \eqref{binipl}, bounds for the initial time $t = 0$, namely
\ber\label{inig}
|G_0[p](0)| &=& \left|\int_0^\infty\int_0^{\play_r[p](0)}\rho(v,r)\dd v\dd r\right|
\le \min\{C_\rho, C_\rho^* \max \{|p(0)|, K\}\}\,,\\ \label{inivg}
|V_0[p](0)| &=& \left|\int_0^\infty\int_0^{\play_r[p](0)}v\rho(v,r)\dd v\dd r\right|
\le C_\rho \max \{|p(0)|, K\}\,.
\eer
The Preisach operator admits also a family of ``nonlinear'' energies.
As a consequence of \eqref{enerpl},
we have for a.~e. $t$ the inequality
\be{epl1}
\play_r[p]_t (p - \play_r[p])\ge 0\,,
\ee
hence
\be{epl2}
\play_r[p]_t (h(p) - h(\play_r[p]))\ge 0
\ee
for every nondecreasing function $h: \real \to \real$. Hence,
for every absolutely continuous input $p$, a counterpart of \eqref{enerpr} in the form
\be{eprh1}
G_0[p]_t h(p) - V_h[p]_t \ge 0\ \mbox{ a.~e.}
\ee
holds with a modified potential
\be{hpot}
V_h[p](t) = \int_0^\infty \int_0^{\play_r[p](t)} h(v) \rho (r,v)\dd v\dd r\,.
\ee
This is related to the fact that for every absolutely continuous nondecreasing function
$\hat h:\real\to \real$, the mapping
$G_{\hat h} := G_0\circ \hat h$ is also a Preisach operator, see \cite{preis}.


\section{Main results}\label{mai}

We denote
\be{2e1}
X_0 = \{\phi \in W^{1,2}(\Omega; \real^3): \phi\big|_{\partial\Omega} =0\}\,, \quad
X = W^{1,2}(\Omega)\,, \quad X_q = W^{1,q}(\Omega)
\ee
for $q > 1$, and reformulate Problem
\eqref{e23}--\eqref{e25} in variational form for all test functions $\phi \in X_0$,
$\psi \in X$, and $\zeta \in X_{q^*}$ for a suitable $q^* > 2$ as follows:
\ber\label{e23v}
\io (\rho_S u_{tt}\cdot\phi + (\BB\nas u_t {+} P[\nas u]):\nas\phi + (p -\beta\theta) \dive\phi)\dd x
&=& \io g\cdot\phi \dd x,\\ \label{e24v}
\io ((G[p]_t - \dive u_t)\psi + \frac{1}{\rho_L} \mu(p)\nabla p \nabla\psi)\dd x
&=& \ipo \gamma_p(x)(p^* - p)\psi\dd s(x),\qquad\\ \nonumber
\io \big((c_0 \theta_t - \|D_P[\nas u]_t\|_* - |D_G[p]_t|)\zeta
 + \kappa(\theta) \nabla\theta\cdot\nabla\zeta\big) \dd x &&\\ \label{e25v}
- \io(\BB\nas u_t:\nas u_t + \frac{1}{\rho_L} \mu(p)|\nabla p|^2-\beta \theta \dive u_t)\zeta\dd x
&=& \ipo \gamma_\theta(x)(\theta^* - \theta)\zeta\dd s(x)\,.
\eer

\begin{hypothesis}\label{h1}
We assume that $\Omega$ is a bounded domain with $C^{1,1}$ boundary. We fix an arbitrary
final time $T>0$, a constant $\bar\theta > 0$, 
and functions $p^* \in W^{1,\infty}(\partial\Omega\times (0,T))$,
$\theta^* \in L^{\infty}(\partial\Omega\times (0,T))$ such that $\theta^*(x,t) \ge \bar\theta$,
$g \in L^2(\Omega\times (0,T))$,
$\gamma_\theta \in L^\infty(\partial\Omega)$,
$\gamma_p \in W^{1,\infty}(\partial\Omega)$, $\gamma_\theta \ge 0$, $\gamma_p \ge 0$ a.~e.,
$\int_0^T\ipo \gamma_p(x)\dd s(x) > 0$, $\beta\in\real$, $c_0>0$.
The coefficients $\rho_S, \rho_L$ are constant and positive, and $\BB$
is the isotropic symmetric positive definite fourth order tensor of the form
\be{hy1}
\BB_{ijkl} = 2\eta \delta_{ik}\delta_{jl} + \omega \delta_{kl}\delta_{ij}
\ee
with constants $\eta>0$, $\omega > 0$.
The nonlinearities in \eqref{e23v}--\eqref{e25v} satisfy the following conditions
\begin{itemize}
\item[{\rm (i)}] $\mu :\real \to [\mu_0, \mu_1]$ is a $C^1$ function,
$0 < \mu_0 < \mu_1$ are fixed constants, and we set
\be{mum}
M(p) = \int_0^p \mu(p')\dd p'.
\ee
\item[{\rm (ii)}] $\kappa :\real \to (0,\infty)$ is a $C^1$ function, $\kappa(0) > 0$,
and there exist constants $0 < a < b < \frac{27}{5} + \frac{12}{5}a$ such that
$$
\liminf_{\theta\to\infty} \frac{\kappa(\theta)}{\theta^{1+a}} > 0\,, \quad
\limsup_{\theta\to\infty} \frac{\kappa(\theta)}{\theta^{1+b}} < \infty.
$$
\item[{\rm (iii)}] $G[p] = f(p) + G_0[p]$, where $G_0$ is the Preisach operator
from Subsection \ref{prei} with an initial memory state $\lambda_{-1} \in \Lambda_K$
for some $K\ge \sup |p^*|$. The dissipation operator $D_G$ associated with $G$
is defined in \eqref{prpot}--\eqref{VG},
and $f: \real \to (0, f_1)$ for some $f_1>0$ is a $C^1$ function such that there exist
$0 < f_2 < f_3$ with the property
\be{hyf}
f_2 \le f'(p)(1+p^2) \le f_3\quad \forall p\in \real\,.
\ee
\item[{\rm (iv)}] Let $\tens$ denote the space of symmetric $3\times 3$ tensors. We assume
that the operator $P:C([0,T];\tens) \to C([0,T];\tens)$ has the form \eqref{P} with $P_0$
defined in Proposition \ref{p3}, and with dissipation operator $D_P$ defined in \eqref{VP}.
\end{itemize}
We prescribe initial conditions \eqref{ini} with
$u^0 \in X_0\cap W^{2,2}(\Omega; \real^3)$, $u^1 \in X_0$,
$p^0 \in L^\infty(\Omega)\cap W^{1,2}(\Omega)$, $|p^0(x)| \le K$ a.~e.,
$\theta^0 \in L^\infty(\Omega)$, $\theta^0(x) \ge \bar\theta$ a.~e.
\end{hypothesis}

Condition (ii) in Hypothesis \ref{h1} is a slight generalization of Hypothesis (I)
of \cite{rr}. We will see the role that it plays in the existence proof.

The main result of this paper reads as follows.

\begin{theorem}\label{t1}
Let Hypotheses \ref{h2}, \ref{h1} hold. Then there exists $q^*$ depending on $a$ and $b$
and a solution $(u,p,\theta)$ to \eqref{e23v}--\eqref{e25v}, and \eqref{ini} with the properties
$u_t \in L^2(0,T; X_0\cap W^{2,2}(\Omega;\real^3))$, $u_{tt} \in L^2(\Omega \times (0,T))$,
$p \in L^\infty(\Omega \times (0,T))$, $M(p) \in L^2(0,T; W^{2,2}(\Omega))$ with
$M(p)$ given by \eqref{mum}, $p_{t} \in L^2(\Omega \times (0,T))$,
$\theta \in L^z(\Omega \times (0,T))$ for every $z < 8 + 3a$,
$\nabla\theta \in L^2(\Omega \times (0,T))$, $\theta_t \in L^2(0,T; W^{-1,q^*}(\Omega))$
for some $q^* > 2$. 
\end{theorem}

We first regularize the problem, prove the existence of a solution for the regularized
system, derive estimates independent of the regularization parameters, and pass to the limit.


\section{Regularization}\label{reg}

We choose regularizing parameters $R >K$ with $K$ from Hypothesis \ref{h2}
and $\delta > 0$ with the intention to let $R \to \infty$ and
$\delta \to 0$, and define mappings $Q_R:\real \to [0,R]$ and $K_R:\real \to \real$ by the formulas
\be{qrkr}
Q_R(z) = \max\{0, \min\{z, R\}\}\,,\quad K_R(z) = \max\{z-R, \min\{0, z+R\}\}\ \for z \in \real\,.
\ee
Let $\Bcal: W^{2,2}(\Omega;\real^3)\cap X_0 \to L^2(\Omega;\real^3)$ denote the mapping
\be{bc}
\Bcal v = -\dive\BB\nas v
\ee 
It follows from a vector counterpart of \cite[Lemma 9.17]{gt},
cf.~also the methods proposed in \cite[Lemma~3.2, p.~260]{necas}, that $|\Bcal v|_2$ is an equivalent norm
for $v$ in $W^{2,2}(\Omega) \cap X_0$, that is, there exist positive constants $C_1 < C_2$
such that for every $v \in W^{2,2}(\Omega;\real^3)\cap X_0$ we have
\be{embe}
C_1 \|v\|_{W^{2,2}(\Omega)} \le |\Bcal v|_2 \le C_2 \|v\|_{W^{2,2}(\Omega)}\,.
\ee
We replace \eqref{e23v}--\eqref{e25v} by the system
\ber \nonumber
\io (\rho_S u_{tt}\cdot\phi + \delta\Bcal u_{tt}\cdot \Bcal\phi + (\BB\nas u_t {+} P[\nas u]):\nas\phi)\dd x &&
\\ \label{e23r} + \io (p -\beta Q_R(\theta)) \dive\phi \dd x
&=& \io g\cdot\phi \dd x,\\ \label{e24r}
\io (((K_R(p)+G[p])_t - \dive u_t)\psi + \frac{1}{\rho_L} \mu(p)\nabla p \nabla\psi)\dd x
&=& \ipo \gamma_p(x)(p^* - p)\psi\dd s(x),\qquad\\ \nonumber
\io \big((c_0 \theta_t - \|D_P[\nas u]_t\|_* - |D_G[p]_t|)\zeta
 + \kappa(\theta) \nabla\theta\cdot\nabla\zeta\big) \dd x &&\\ \label{e25r}
- \io(\BB\nas u_t{:}\nas u_t {+} \frac{1}{\rho_L} \mu(p)Q_R(|\nabla p|^2) {-}\beta Q_R(\theta) \dive u_t)\zeta\dd x
&=& \ipo \gamma_\theta(x)(\theta^* - \theta)\zeta\dd s(x)
\eer
with test functions $\phi \in W^{2,2}(\Omega;\real^3)\cap X_0$, $\psi, \zeta \in X$
and initial conditions \eqref{ini}.

\begin{proposition}\label{p4}
In addition to Hypotheses \ref{h2}, \ref{h1}, assume that $u^1 \in X_0\cap W^{2,2}(\Omega; \real^3)$.
Then there exists a solution $(u,p,\theta)$ to \eqref{e23r}--\eqref{e25r} with the properties
$u_{tt} \in L^2(0,T; W^{2,2}(\Omega;\real^3)\cap X_0)$,
$p \in L^2(0,T; W^{1,2}(\Omega))$, $M(p) \in L^2(0,T; W^{2,2}(\Omega))$ with
$M(p)$ given by \eqref{mum}, $p_{t} \in L^2(\Omega \times (0,T))$,
$\theta \in L^2(\Omega \times (0,T))$,
$\nabla\theta \in L^2(\Omega \times (0,T))$, $\theta_t \in L^2(0,T; X^*)$, where $X^*$ is the dual of $X$.
\end{proposition}

System \eqref{e23r}--\eqref{e25r} for each fixed $R>0$ and $\delta > 0$
will be solved by Faedo-Galerkin approximations.
We choose $\Ecal=\{e_k; k=1,2,\dots\}$ in $L^2(\Omega;\real^3)$ and
$\Wcal=\{w_j; j=0,1,2,\dots\}$ in $L^2(\Omega)$
to be the complete orthonormal systems of eigenfunctions defined by
\be{eigen}
\Bcal e_k = \lambda_k e_k \ \text{ in } \ \Omega\,, \ \ e_k\big|_{\partial\Omega} = 0\,,
\quad -\Delta w_j = \mu_j w_j \ \text{ in } \ \Omega\,, \ \ \nabla w_j\cdot n\big|_{\partial\Omega} = 0\,,
\ee
with $\mu_0 = 0$, $\lambda_k>0, \mu_j > 0$ for $j,k\ge 1$, and put for $n\in\nat$
\be{faeg}
u\om(x,t) = \sum_{k=1}^n u_k(t) e_k(x)\,, \quad \theta\om(x,t) = \sum_{j=0}^n \theta_j(t) w_j(x)
\ee
with coefficients $u_k:[0,T] \to \real$, $\theta_j:[0,T] \to \real$ which will be determined
as the solution of the system
\ber \nonumber
(\rho_S {+} \lambda_k {+} \delta \lambda_k^2)\ddot u_k + \io P[\nas u\om]:\nas e_k \dd x &&
\\ \label{e23fgk} + \io (p\om -\beta Q_R(\theta\om)) \dive e_k\dd x
&=& \io g\cdot e_k \dd x,\\ \nonumber
\io ((K_R(p\om)+G[p\om])_t - \dive u\om_t)\psi \dd x &&\\ \label{e24fgk}
+ \frac{1}{\rho_L} \io \mu(p\om)\nabla p\om\cdot \nabla\psi\dd x
&=& \ipo \gamma_p(x)(p^* - p\om)\psi\dd s(x),\qquad\\ \nonumber
c_0 \dot \theta_j + \io \big({-} |D_G[p\om]_t| w_j
 {+} \kappa(\theta\om) \nabla\theta\om\cdot\nabla w_j\big) \dd x &&\\ \nonumber
+ \io(\beta Q_R(\theta\om) \dive u\om_t - \|D_P[\nas u\om]_t\|_*)w_j \dd x &&\\
\label{e25fgk}
- \io(\BB\nas u\om_t{:}\nas u\om_t {+} \frac{1}{\rho_L} \mu(p\om)Q_R(|\nabla p\om|^2)) w_j \dd x
&=& \ipo \gamma_\theta(x)(\theta^* - \theta\om) w_j \dd s(x)
\eer
for $k=1, \dots, n$ and $j=0,1,\dots, n$, and for all $\psi \in X$.
We prescribe initial conditions
\ber\label{inidu}
u_k(0) = \io u^0(x)\cdot e_k(x) \dd x\,, && \dot u_k(0) = \io u^1(x)\cdot e_k(x) \dd x\,,\\ \label{inidp}
\theta_j(0) = \io \theta^0(x) w_j(x) \dd x\,, && p\om(x,0) = p^0(x)\,.
\eer
This is an ODE system \eqref{e23fgk}, \eqref{e25fgk} coupled with a standard PDE
with hysteresis \eqref{e24fgk},
which has a strong solution in a maximal interval of existence, which coincides with the whole
interval $[0,T]$ provided we prove that the solution remains bounded in the maximal
interval of existence.

Put $\Ecal_n=\{e_k; k=1,2,\dots, n\}$ and $\Wcal_n=\{w_j; j=0,1,2,\dots, n\}$.
Then \eqref{e23fgk}--\eqref{e25fgk} can be equivalently written as
\ber \nonumber
\io (\rho_S u\om_{tt}\cdot\phi + (\BB\nas u\om_t {+} P[\nas u\om]):\nas\phi)\dd x &&
\\ \label{e23fg} + \delta\io \Bcal u\om_{tt}\cdot \Bcal\phi \dd x
+ \io (p\om -\beta Q_R(\theta\om)) \dive\phi \dd x
&=& \io g\cdot\phi \dd x,\\ \nonumber
\io ((K_R(p\om)+G[p\om])_t - \dive u\om_t)\psi \dd x &&\\ \label{e24fg}
 + \io \frac{1}{\rho_L} \mu(p\om)\nabla p\om \nabla\psi\dd x
&=& \ipo \gamma_p(x)(p^* - p\om)\psi\dd s(x),\qquad\\ \nonumber
\io \big((c_0 \theta\om_t - |D_G[p\om]_t|)\zeta
 + \kappa(\theta\om) \nabla\theta\om\cdot\nabla\zeta\big) \dd x &&\\ \nonumber
 + \io(\beta Q_R(\theta\om)\dive u\om_t - \|D_P[\nas u\om]_t\|_*)\zeta\dd x&& \\ \label{e25fg}
- \io(\BB\nas u\om_t{:}\nas u\om_t + \frac{1}{\rho_L} \mu(p\om)Q_R(|\nabla p\om|^2))\zeta\dd x
&=& \ipo \gamma_\theta(x)(\theta^* - \theta\om)\zeta\dd s(x)
\eer
with test functions $\phi \in \mathrm{Span}~\Ecal_n$, $\zeta \in \mathrm{Span}~\Wcal_n$, and $\psi\in X$.

We now derive a series of estimates. By $C$ we denote any positive constant depending only on the data,
by $C_R$ any constant depending on the data and on $R$, and by $C_{R,\delta}$ any constant
depending on the data, on $R$, and on $\delta$, all independent of the dimension $n$ of the Galerkin
approximation.

To simplify the presentation, we introduce from now on the notation $|\cdot|_q$ for the norm
in $L^q(\Omega)$, and by $\|\cdot\|_q$ the norm in $L^q(\Omega\times (0,T))$. We will systematically
use the Gagliardo-Nirenberg inequality in the form
\be{gn}
|w|_q \le C(|w|_s + |w|_s^{1-\gamma} |\nabla w|_r^\gamma)\,, \quad \gamma = \frac{\frac{1}{s} - \frac{1}{q}}
{\frac{1}{3} + \frac{1}{s} - \frac{1}{r}}
\ee
which holds for every $w \in W^{1,r}(\Omega)$ and every $\frac{1}{s} > \frac{1}{q} > \frac{1}{r} - \frac{1}{3}$.
For the proof, see, e.\,g., \cite[\S 15]{bin}.


\subsection{Estimate 1}\label{estm1}

We test \eqref{e23fg} with $\phi=u_t^{(n)}$ and \eqref{e24fg} with $\psi = p\om$ and sum up the results
to obtain
\ber \nonumber
&&\hspace{-16mm}\io \Big(\rho_S u\om_{tt}\cdot u\om_t + \delta \Bcal u\om_{tt}\cdot \Bcal u\om_t
+ (\BB\nas u\om_t + P[\nas u\om]):\nas u\om_t \\ \nonumber
&& +\, (K_R(p\om)+G[p\om])_t p\om + \frac{1}{\rho_L} \mu(p\om) |\nabla p\om|^2\Big)\dd x\\ \label{em7}
&=& \io (\beta Q_R(\theta\om) \dive u\om_t + g\cdot u\om_t)\dd x + \ipo\gamma_p(x)p\om (p^* - p\om)\dd s(x).
\eer
Integrating in time from $0$ to $t$ and using the energy identities \eqref{e10} we obtain for all $t\in (0,T)$
the estimate
\ber \nonumber
&&\hspace{-16mm}|u\om_t(t)|_2^2 + \delta |\Bcal u\om_t(t)|_2^2 + |\nas u\om(t)|_2^2
+ \|\nas u\om_t\|_2^2\\ \label{em1}
&& +\, |p\om(t)|_2^2 + \|\nabla p\om\|_2^2 + \int_0^T\ipo \gamma_p(x) |p\om|^2\dd s(x)\dd t \le C_R\,.
\eer


\subsection{Estimate 2}\label{estm2}

We choose in \eqref{e24fg} $\psi = M(p\om)_t$ with $M$ given by \eqref{mum}.
By Proposition \ref{p5} and formula \eqref{pre}, we have
$G_0[p\om]_t M(p\om)_t = G_0[p\om]_t p\om_t \mu(p\om) \ge 0$.
By Hypothesis \ref{h1}\,(i), (iii), we thus have
$$
(K_R(p\om)+ G[p\om]_t) M(p\om)_t
\ge \mu_0 \min\left\{1, \frac{f_2}{1+R^2}\right\} |p\om_t|^2
$$
and we obtain for all $t\in (0,T)$ the estimate
\be{em4}
\|p\om_t\|_2^2 + |\nabla p\om(t)|_2^2 + \ipo \gamma_p(x) |p\om(t)|^2\dd s(x) \le C_R\,.
\ee
By comparison in Eq.~\eqref{e24fg}, we see that
\be{em4a}
\|M(p\om)\|_{L^2(0,T; W^{2,2}(\Omega))} \le C_R\,.
\ee


\subsection{Estimate 3}\label{estm3}

Choosing in \eqref{e23fg} $\phi=u_{tt}^{(n)}$ yields
\ber \nonumber
&&\hspace{-16mm}\io \Big(\rho_S |u\om_{tt}|^2 + \delta |\Bcal u\om_{tt}|^2
+ \BB\nas u\om_t:\nas u\om_{tt} + \frac{\partial}{\partial t}(P[\nas u\om]:\nas u\om_t)\Big)\dd x
\\ \label{em2}
&=& \io \Big(P[\nas u\om]_t:\nas u\om_t +(\beta Q_R(\theta\om)-p^{(n)}) \dive u\om_{tt} + g\cdot u\om_{tt}\Big)\dd x.
\eer
We now integrate in time again and use Proposition \ref{p3}, estimate \eqref{em4},
as well as the Gronwall argument, to conclude for all $t\in (0,T)$ that
\be{em3} 
\|u\om_{tt}\|_2^2 + \delta \|\Bcal u\om_{tt}\|_2^2 + |\nas u\om_t(t)|_2^2 \le C_{R,\delta}\,.
\ee


\subsection{Estimate 4}\label{estm4}

We choose in \eqref{e25fg} $\zeta = \theta^{(n)}$. The only superlinear term in \eqref{e25fg}
is $\BB\nas u\om_t:\nas u\om_t$ which has to be estimated in the norm of $L^2(\Omega\times (0,T))$,
that is, $\nas u\om_t$ has to be estimated in $L^4(\Omega\times (0,T))$. This will be done
using the Gagliardo-Nirenberg inequality \eqref{gn}, which yields for every $t \in (0,T)$ that
\be{em5}
|\nas u\om_t(t)|_4 \le C(|\nas u\om_t(t)|_2 + |\nas u\om_t(t)|_2^{1/4} |\Bcal u\om_t(t)|_2^{3/4}) \le C_{R,\delta}
\ee
by virtue of \eqref{em1} and \eqref{embe}. Note also that we have the pointwise inequalities
$$
|D_G[p\om]_t| \le C |p\om_t|\,, \quad \|D_P[\nas u\om]_t\|_* \le C |\nas u\om_t|
$$
which follow from \eqref{prpot}, \eqref{mono}, \eqref{VP}, \eqref{h4b}, and Proposition \ref{p3}\,(i).
We thus obtain
\be{em6}
|\theta\om(t)|_2^2 + \|\nabla\theta\om\|_2^2 + \int_0^T \ipo \gamma_\theta(x) |\theta\om|^2 \dd s(x)\dd t \le C_{R,\delta}
\ee
for all $t \in (0,T)$. Finally, let $\zeta \in L^2(0,T; X)$ be arbitrary,
$\zeta(x,t) = \sum_{j=0}^\infty \zeta_j(t)w_j(x)$.
We test \eqref{e25fg} with $\zeta = \zeta_j(t)$ and obtain using the previous estimates that
\be{em8}
\int_0^T\io \theta\om_t\zeta\dd x\dd t \le C_{R,\delta}\|\zeta\|_{L^2(0,T; X)}\,,
\ee
or, in other words,
\be{em9}
\|\theta\om_t\|_{L^2(0,T; X^*)} \le C_{R,\delta}\,.
\ee


\subsection{Passage to the limit as $n \to \infty$}\label{ninf}

We keep for the moment the regularization parameters $R$ and $\delta$ fixed, and let $n$ tend to $\infty$.
By a standard argument based on compact anisotropic embeddings, see \cite{bin}, we infer,
passing to a subsequence, if necessary, that there exist functions $(u,p,\theta)$ such
that the following convergences take place:
$$
\begin{array}{rclrl}
u\om_{tt} & \to & u_{tt} & \mbox{ weakly in } & L^2(0,T; W^{2,2}(\Omega;\real^3)\cap X_0)\,,\\
\nas u\om_t & \to & \nas u_t & \mbox{ strongly in } & L^4(\Omega; C([0,T];\tens))\,,\\
\nas u\om & \to & \nas u & \mbox{ strongly in } & L^4(\Omega; C([0,T];\tens))\,,\\
P[\nas u\om] & \to & P[\nas u] & \mbox{ strongly in } & L^4(\Omega; C([0,T];\tens))\,,\\
\|D_P[\nas u\om]_t\|_* & \to & \|D_P[\nas u]_t\|_* & \mbox{ strongly in } & L^2(\Omega; C[0,T])\,,\\
p\om & \to & p & \mbox{ strongly in } & L^4(\Omega; C[0,T])\,,\\
p\om_t & \to & p_t & \mbox{ weakly in } & L^2(\Omega\times (0,T))\,,\\
K_R(p\om)_t & \to & K_R(p)_t & \mbox{ weakly in } & L^2(\Omega\times (0,T))\,,\\
G[p\om]_t & \to & G[p]_t & \mbox{ weakly in } & L^2(\Omega\times (0,T))\,,\\
|D_G[p\om]_t| & \to & |D_G[p]_t| & \mbox{ weakly in } & L^2(\Omega\times (0,T))\,,\\
\nabla p\om & \to & \nabla p & \mbox{ strongly in } & L^2(\Omega\times (0,T);\real^3)\,,\\
Q_R(|\nabla p\om|^2) & \to & Q_R(|\nabla p|^2) & \mbox{ strongly in } & L^2(\Omega\times (0,T))\,,\\
\gamma_p p\om & \to & \gamma_p p & \mbox{ strongly in } & L^2(\partial\Omega \times [0,T])\,,\\
\theta\om & \to & \theta & \mbox{ strongly in } & L^2(\Omega\times (0,T))\,,\\
\theta\om_t & \to & \theta_t & \mbox{ weakly in } & L^2(0,T; X^*)\,,\\
\nabla\theta\om & \to & \nabla\theta & \mbox{ weakly in } & L^2(\Omega\times (0,T);\real^3)\,,\\
\gamma_\theta \theta\om & \to & \gamma_\theta \theta & \mbox{ strongly in } & L^2(\partial\Omega \times [0,T])\,.
\end{array}
$$
The convergences of the hysteresis terms $P[\nas u\om]$, $G[p\om]_t$,
$\|D_P[\nas u\om]_t\|_*$, $|D_G[p\om]_t|$ follow indeed from
\eqref{P}, \eqref{h8}, \eqref{G}, \eqref{lipg},
and \eqref{e10}. We can therefore let $n$ tend to $\infty$ in \eqref{e23r}--\eqref{e25r} and conclude
that the limit $(u, p, \theta)$ satisfies the conditions of Proposition~\ref{p4}.


\section{Proof of Theorem \ref{t1}}\label{pt1}

In this section, we show that a sequence of solutions to \eqref{e23r}--\eqref{e25r} converges to a solution
to \eqref{e23v}--\eqref{e25v} as $R\to \infty$ and $\delta \to 0$. To this end, we fix sequences
$\{R_i\}, \{\delta_i\}$ for $i \in \nat$ such that
\be{ep1}
\lim_{i\to \infty} R_i = \infty\,, \quad \lim_{i\to \infty} \delta_i = 0\,, 
\ee
and choose a sequence $\{u^1_i\}$ of initial conditions in $X_0\cap W^{3,2}(\Omega; \real^3)$ such that
\be{ep2}
\lim_{i\to \infty} \|u^1_i - u^1\|_{X_0} = 0\,,
\quad \lim_{i\to \infty} \delta_i \|u^1_i\|^2_{W^{3,2}(\Omega; \real^3)} = 0\,.
\ee
We further denote by $(u\oi, p\oi, \theta\oi)$ solutions $(u, p, \theta)$ to Problem \eqref{e23r}--\eqref{e25r}
corresponding to the choice $R = R_i$, $\delta = \delta_i$, $u^1 = u^1_i$. The next step
is to derive some properties of the sequence $(u\oi, p\oi, \theta\oi)$ independent of $i$.


\subsection{Positivity of temperature}\label{posi}

We first observe that there exists a constant $C>0$ such that for every nonnegative test function $\zeta \in X$
we have by virtue of \eqref{e25r} that
\be{e25p}
\io (c_0 \theta\oi_t \zeta
 + \kappa(\theta) \nabla\theta\oi\cdot\nabla\zeta) \dd x \ge 
- C\io Q_{R_i}(\theta\oi)^2\zeta\dd x
+ \ipo \gamma_\theta(x)(\theta^* - \theta\oi)\zeta\dd s(x).
\ee
Let $v(t)$ be the solution of the ODE
\be{eqv}
c_0 \dot v(t) = -Cv^2(t)\,, \quad v(0) = \bar\theta\,,
\ee
that is,
\be{eqv2}
v(t) = \left(\frac{C}{c_0}t + \frac{1}{\bar\theta}\right)^{-1}.
\ee
For every nonnegative test function $\zeta \in X$ we have in particular
\be{e25p2}
\io (c_0 v_t \zeta
 + \kappa(\theta\oi) \nabla v\cdot\nabla\zeta) \dd x \le 
- C\io v^2\zeta\dd x
+ \ipo \gamma_\theta(x)(\theta^* - v)\zeta\dd s(x).
\ee
Subtracting \eqref{e25p} from \eqref{e25p2} we obtain
\be{e25p1}
\io (c_0 (v - \theta\oi)_t \zeta
 + \kappa(\theta\oi) \nabla (v - \theta\oi)\cdot\nabla\zeta) \dd x \le 
 C\io (Q_{R_i}^2(\theta\oi) - v^2)\zeta\dd x
+ \ipo \gamma_\theta(x)(\theta\oi - v)\zeta\dd s(x).
\ee
We now choose any smooth convex function $F: \real \to \real$ such that $F(s) = 0$ for $s\le 0$,
$F(s) > 0 $ for $s>0$, and test \eqref{e25p1} by $\zeta = F'(v-\theta\oi)$. We have in all cases
$$
(Q_{R_i}^2(\theta\oi) - v^2)F'(v-\theta\oi) \le 0 \ \ale,
$$
hence
\be{e25p3}
c_0 \frac{\dd}{\dd t}\io F(v - \theta\oi) \dd x \le 0\,,
\ee
and we conclude for every $i \in \nat$ that
\be{ee0}
\theta\oi(x,t) \ge v(t) \ \ale 
\ee
We now pass to a series of estimates independent of $i$.
To simplify the presentation, we occasionally omit the indices $\{\}\oi$
in the computations in Subsections \ref{est1}--\ref{est7} below,
and write simply $(u, p, \theta)$ instead of $(u\oi, p\oi, \theta\oi)$
whenever there is no risk of confusion.
The symbol $C$ denotes as before any constant independent of $i$.


\subsection{Estimate 5}\label{est1}

Test \eqref{e23r} by $\phi = u_t = u\oi_t$, \eqref{e24r} by $\psi = p = p\oi$, and \eqref{e25r} by $\zeta = 1$.
Summing up the three resulting equations we obtain by virtue of \eqref{e10} that
\ber\nonumber
&&\hspace{-14mm}\frac{\dd}{\dd t} \io\left(c_0\theta + \frac{\rho_S}{2} |u_t|^2
+ \frac{\delta_i}{2} |\Bcal u_t|^2 + V_P[\nas u] + V_G[p] + \hat K_{R_i}(p)\right)\dd x\\ \nonumber
&&+\, \frac{1}{\rho_L} \io \mu(p)\big(|\nabla p|^2 {-} Q_{R_i}(|\nabla p|^2)\big)\dd x
\\ \label{ee1}
&=& \io g\cdot u_t\dd x + \ipo \left(\gamma_p(x)(p^* - p) p
+ \gamma_\theta(x)(\theta^* - \theta)\right) \dd s(x)\,,
\eer
where we set $\hat K_R(p) = \int_0^p K_R'(p')p'\dd p'$ for $p\in\real$ and $R_i>0$. Integrating in time
and using \eqref{VP} and \eqref{ep2}, we get for every $t \in (0,T)$ the estimate
\be{ee2}
\io\left(\theta\oi + |u\oi_t|^2 + |\nas u\oi|^2 \right)(x,t)\dd x
+ \int_0^T\!\ipo \left(\gamma_\theta(x)\theta\oi + \gamma_p(x)|p\oi|^2\right)(x,t) \dd s(x)\dd t \le C\,.
\ee


\subsection{Estimate 6}\label{est2}

We test \eqref{e25r} by $\zeta = \theta^{-a} = (\theta\oi)^{-a}$
with $a$ from Hypothesis \ref{h1}\,(ii), and observe that
\ber \nonumber
&&\hspace{-15mm}\io \theta^{-a}\left(\|D_P[\nas u]_t\|_* + |D_G[p]_t| + \BB\nas u_t{:}\nas u_t
+ \frac{1}{\rho_L} \mu(p)Q_{R_i}(|\nabla p|^2)\right)\dd x\\ \nonumber
&& +\, a \io \kappa(\theta)\theta^{-1-a} |\nabla\theta|^2\dd x
\\ \label{e25a}
&=& \beta\io Q_{R_i}(\theta)\theta^{-a} \dive u_t\dd x
- \ipo \gamma_\theta(x) \theta^{-a}(\theta^* - \theta)\dd s(x)
+ \frac{c_0}{1-a}\frac{\dd}{\dd t} \io \theta^{1-a}\dd x.
\eer
Integrating in time and using \eqref{ee0}, \eqref{ee2}, and Hypothesis \ref{h1}\,(ii), we obtain in particular
\be{ee3}
\int_0^T\io (\theta^{-a}|\dive u_t|^2 + |\nabla\theta|^2)\dd x\dd t
\le C\left(1+ \int_0^T\io \theta^{1-a} |\dive u_t|\dd x\dd t\right).
\ee
The integral on the right hand side can be estimated by H\"older's inequality
$$
\int_0^T\io \theta^{1-a} |\dive u_t|\dd x\dd t \le \left(\int_0^T\io \theta^{-a} |\dive u_t|^2\dd x\dd t\right)^{1/2}
\left(\int_0^T\io \theta^{2-a} \dd x\dd t\right)^{1/2},
$$
which entails that
\be{ee4}
\int_0^T\io |\nabla\theta|^2\dd x\dd t \le C\left(1+ \int_0^T\io \theta^{2-a}\dd x\dd t\right).
\ee
Applying the Gagliardo-Nirenberg inequality \eqref{gn} with $s=1$, $r=2$, and $q=2-a$
and using \eqref{ee2}, we estimate the right hand side of \eqref{ee4} from above by
$C(1+\|\nabla\theta\|_2^{(1-a)6/5})$. Thus, for $H :=\|\nabla\theta\|_2^2$, inequality \eqref{ee4}
has the form $H \le C (H^\omega + 1)$ with $\omega = (1-a)3/5 < 1$, and
$H^\omega \le \omega\delta H + (1-\omega)\delta^{-\omega/(1-\omega)}$. Choosing, for example,
$\delta=1/C$, we obtain
\be{ee5}
\|\nabla\theta\oi\|_2 \le C\,.
\ee
Using \eqref{gn} again with $s=1$, $r=2$, and $q=8/3$, we obtain
\be{ee6}
\|\theta\oi\|_{8/3} \le C\,.
\ee


\subsection{Estimate 7}\label{est3}

Test \eqref{e23r} by $\phi = u_t = u\oi_t$ and \eqref{e24r} by $\psi = p = p\oi$.
The sum of the two equations yields
\ber \nonumber
&&\hspace{-10mm}\frac{\dd}{\dd t}\io \left(\frac{\rho_S}{2} |u_t|^2 + \frac{\delta_i}{2}|\Bcal u_t|^2
+ \hat K_R(p)\right)\dd x \\ \nonumber
&&+\, \io\left((\BB\nas u_t + P[\nas u]):\nas u_t
+ G[p]_t p + \frac{1}{\rho_L} \mu(p) |\nabla p|^2\right)\dd x\\ \label{ee7}
&& = \io (\beta Q_{R_i}(\theta)\dive u_t + g\cdot u_t)\dd x + \ipo\gamma_p(x)p (p^* - p)\dd s(x).
\eer
By \eqref{e10} and \eqref{VG}, we have $G[p]_t p \ge V_G[p]_t$, $V_G[p](x,t) \ge 0$, and
$V_G[p](x,0) \le C|p(x,0)|^2 \le C$ a.~e.
Integrating in time, taking into account \eqref{ep2} and the previous estimates, we get
\be{ee8}
\|\nas u\oi_t\|_2 + \|\nabla p\oi\|_2 + \delta_i|\Bcal u\oi_t(t)|_2
+ \int_0^T\ipo \gamma_p(x)|p\oi|^2\dd s(x) \le C
\ee
for all $t \in (0,T)$. Consequently, as $\gamma_p$ does not identically vanish
on $\partial\Omega$ by Hypothesis \ref{h1}, we also have
\be{ee8a}
\|p\oi\|_{L^2(0,T;W^{1,2}(\Omega))} \le C.
\ee


\subsection{Estimate 8}\label{est4}

Test \eqref{e23r} by $\phi = u_{tt} = u\oi_{tt}$. Then
\ber \nonumber
&&\hspace{-10mm}\io (\rho_S |u_{tt}|^2 + \delta_i |\Bcal u_{tt}|^2)\dd x
+ \frac12 \frac{\dd}{\dd t}\io(\BB\nas u_t {+} 2P[\nas u]):\nas u_t \dd x\\ \label{ee9}
&&\le C\io(|g| + |\nabla p| + |\nabla\theta|)|u_{tt}|\dd x + \io P[\nas u]_t:\nas u_t \dd x.
\eer
{}From Proposition \ref{p3}, \eqref{ee5}, Korn's inequality, and \eqref{ee8} it follows that
\be{ee10}
\|u\oi_{tt}\|_2^2 + \delta_i \|\Bcal u\oi_{tt}\|_2^2 + |\nas u\oi_t(t)|_2^2 \le C
\ee
for every $t \in (0,T)$.


\subsection{Estimate 9}\label{est9}

We rewrite \eqref{e23r} in the form
\be{ep3}
\io \big((\rho_S u_{tt}+ \Bcal u_t)\cdot\phi + \delta_i\Bcal u_{tt}\cdot \Bcal\phi\big)\dd x
= \io (f+h)\cdot\phi \dd x
\ee
for all $\phi \in W^{2,2}(\Omega;\real^3)\cap X_0$, where
\be{ep4}
f = g-\beta\nabla Q_R(\theta) + \nabla p\,, \quad h = \dive P[\nas u]\,.
\ee
We have $f \in L^2(\Omega;\real^3)$ by \eqref{ee5}, \eqref{ee8}, and Hypothesis \ref{h1}.
To estimate $h$ in $L^2$, we use
\eqref{h8} and proceed as follows. Let $E_l$, $l=1,2,3$, be the $l$-th coordinate vector,
let $(x,t) \in \Omega\times (0,T)$ be an arbitrary Lebesgue point of $\partial_{x_l}P[\nas u]$,
and let $s_0\in \real$ be sufficiently small such that
$x + sE_l \in \Omega$ for $|s| < s_0$. By \eqref{h8} and \eqref{P} we have
\ber \nonumber
&&\hspace{-15mm}|P[\nas u](x+sE_l) - P[\nas u](x,t)| \le C\Big(|\nas u^0(x+sE_l) - \nas u^0(x)|\\ \label{ep5}
&&+\, \int_0^t |\nas u_t(x+sE_l,\tau) - \nas u_t(x,\tau)|\dd\tau\Big)\,,
\eer
so that in the limit as $s\to 0$ we have
\be{ep5a}
\left|\frac{\partial}{\partial x_l}P[\nas u](x,t)\right|
\le C\left(\left|\frac{\partial}{\partial x_l}\nas u^0(x)\right|
+ \int_0^t \left|\frac{\partial}{\partial x_l}\nas u_t(x,\tau)\right|\dd\tau\right) \ \ale
\ee
and
\be{ep6}
|h(t)|_2 \le C\left(1 + \int_0^t|\Bcal u_t(\tau)|\dd \tau\right).
\ee
Consider now the Fourier expansion of $u = u\oi$ in the form
\be{four}
u(x,t) = \sum_{k=1}^\infty u_k(t) e_k(x)
\ee
similar to \eqref{faeg}, with coefficients
\be{fourc}
u_k(t) = \io u(t)\cdot e_k(x) \dd x\,.
\ee
It follows e.\,g. from \eqref{ee8} that the series
\be{fourb}
\Bcal u_t(x,t) = \sum_{k=1}^\infty \lambda_k \dot u_k(t) e_k(x)
\ee
is strongly convergent in $L^2(\Omega; \real^3)$.

We now test \eqref{ep3} by $\phi = \Bcal u\om_t$, where $u^{(n)}$ is as in \eqref{faeg}
with coefficients $u_k(t)$ given by \eqref{fourc}. Then
\ber \nonumber
&&\hspace{-15mm}\frac{\dd}{\dd t} \io \left(\frac{\rho_S}{2}\BB\nas u\om_t:\nas u\om_t + \frac{\delta_i}{2}
\BB\nas \Bcal u\om_t:\nas \Bcal u\om_t\right)\dd x + |\Bcal u\om_t(t)|_2^2\\ \label{ep7}
&\le& C\left(1 + |f(t)|_2 + \int_0^t|\Bcal u_t(\tau)|_2\dd \tau\right)|\Bcal u\om_t(t)|_2\,,
\eer
hence,
\ber \nonumber
&&\hspace{-15mm}\frac{\dd}{\dd t} \io \left(\frac{\rho_S}{2}\BB\nas u\om_t:\nas u\om_t + \frac{\delta_i}{2}
\BB\nas \Bcal u\om_t:\nas \Bcal u\om_t\right)\dd x + |\Bcal u\om_t(t)|_2^2\\ \label{ep7a}
&\le& C\left(1 + |f(t)|^2_2 + \int_0^t|\Bcal u_t(\tau)|^2_2\dd \tau\right)\,.
\eer
By \eqref{ep2}, we can integrate this inequality from $0$ to $t$, pass to the limit as $n \to \infty$, and use 
Gronwall's argument to obtain in particular that
\be{ee11}
\|\Bcal u_t\|_2= \|\Bcal u\oi_t\|_2 \le C\,.
\ee
The next computation based on \eqref{gn} is to check that
\be{ee12}
\|u\oi_t\|_{L^r(0,T; C(\bar\Omega; \real^3))} \le C \quad \mbox{ for every }\ r\in [1,4)\,.
\ee
Indeed, we choose any $\alpha \in [0, 1/6)$ and put $\frac{1}{q} = \frac13 - \alpha$.
By \eqref{ee10}, \eqref{ee11}, and \eqref{gn} we have
\be{ee13}
|\partial_{x_j} u\oi_t(t)|_q \le C\left(|\partial_{x_j} u\oi_t(t)|_2 + |\partial_{x_j} u\oi_t(t)|_2^{1-\gamma}
|\partial_{x_j} \nabla u\oi_t(t)|_2^\gamma\right)\,,\quad \gamma = \frac{\frac12 - \frac{1}{q}}{\frac13}.
\ee
Then $\partial_{x_j} u\oi_t$ is bounded in $L^p(0,T;L^q(\Omega; \real^3))$ for $p\gamma = 2$, that is,
\be{ee14}
\big|\partial_{x_j} u\oi_t\big|_{L^p(0,T;L^q(\Omega; \real^3))} \le C \ \for q = \frac{3}{1 - 3\alpha},
\ p = \frac{4}{6\alpha + 1}\,.
\ee
By \eqref{ee10}, $u\oi_t$ is bounded in $L^\infty(0,T; L^2(\Omega; \real^3))$,
and by \eqref{gn} for $\alpha >0$ we have
\be{ee15}
|u\oi_t(t)|_\infty \le C\left(|u\oi_t(t)|_2 + |u\oi_t(t)|_2^{1-\hat\gamma}
|\nabla u\oi_t(t)|_q^{\hat\gamma}\right)\,,\quad \hat \gamma = \frac{\frac12}{\frac56 - \frac{1}{q}}.
\ee
We then obtain \eqref{ee12} for $r \hat\gamma = p$, that is,
\be{ee16}
r = 4\frac{1+2\alpha}{1+6\alpha} \in [1,4)\,.
\ee


\subsection{Estimate 10}\label{est5}

In this subsection, we prove the following statement.

\begin{proposition}\label{p2}
Let Hypothesis \ref{h1} hold. Then there exists a constant $C^* > 0$ such that
$$
|p\oi(x,t)| \le C^* \quad \mbox{a.~e.}
$$
\end{proposition}

Note that by \eqref{ee12}, we have
\be{m3}
U \in L^3(0,T)\,, \mbox{ where we put }\ U(t) := 1+ \sup_{x\in \Omega} |u\oi_t(x,t)|\,.
\ee 
As a preliminary step before we pass to the proof of Proposition \ref{p2},
we prove the following auxiliary result for $p = p\oi$.

\begin{lemma}\label{l1}
There exist constants $c>0$ and $C>0$ independent of $m$ such that for every $m\ge 1$
and every $t \in [0,T]$ we have
\ber \nonumber
&&\hspace{-16mm}\io|p(x,t)|^{2{m}}\dd x + c \int_0^t \io |\nabla (p|p|^{m-1})|^2 \dd x\dd\tau\\ \label{m6}
&\le& C(1+m^2) \int_0^t U^2(\tau)\left(K^{2m} + \io|p(x,\tau)|^{2{m}}\dd x\right)\dd\tau\,.
\eer
\end{lemma}

\bpf{Proof}\
We choose an arbitrary ${Q}>0$ and $m\ge 1$, and test  \eqref{e24r} by $h_{{Q},m}(p)$, where we set
\be{pl1}
h_{{Q},m}(p) =
\left\{
\begin{array}{ll}
p |p|^{2m} & \mbox{ for } \ |p| < {Q}\,,\\
{Q}^{2m+1} + (2m+1) (p-{Q}) {Q}^{2m} & \mbox{ for } \ p\ge {Q}\,,\\
-{Q}^{2m+1} + (2m+1) (p+{Q}) {Q}^{2m} & \mbox{ for } \ p\le -{Q}\,.
\end{array}
\right.
\ee
We have $h_{{Q},m}(p) \in L^2(0,T;W^{1,2}(\Omega))$ by virtue of \eqref{ee8a},
hence this is an admissible test function.
By \eqref{eprh1} we have $G_0[p]_t h_{{Q},m}(p) \ge V_{h_{{Q},m}}[p]_t$
and $K_R(p)p_t h_{Q,m}(p) = \hat K_{R,Q,m}(p)_t$, with $\hat K_{R,Q,m}(p)
= \int_0^p K_R(p') h_{Q,m}(p')\dd p'$. We thus have
\ber\nonumber
&&\hspace{-16mm}\frac{\dd}{\dd t} \io (\hat K_{R,Q,m}(p)+ V_{h_{{Q},m}}[p])\dd x
+ \io f'(p)h_{{Q},m}(p) p_t\dd x
+ \mu_0 (2{m}+1) \io |\nabla p|^2 \min\{|p|, {Q}\}^{2{m}}\dd x\\ \label{m2}
&&\le  -(2{m}+1)\io u_t \min\{|p|,{Q}\}^{2{m}}\nabla p \dd x + \ipo \gamma(x) (p^* - p)h_{{Q},m}(p)\dd s(x)\,,
\eer
together with $\hat K_{R,Q,m}(p)(x,0) = 0$ by Hypothesis \ref{h2}, and
\be{m2a}
\io V_{h_{{Q},m}}[p](x,0)\dd x \le CK^{2m}
\ee
as a consequence of \eqref{hpot} and \eqref{binipl},
with $C$ independent of ${Q}$ and $m$. We estimate the right hand side of \eqref{m2} as follows
\bers
&&\hspace{-16mm}-(2{m}+1)\io u_t \min\{|p|,{Q}\}^{2m}\nabla p \dd x\\
&\le& (2{m}+1)U(t)\left(\io \min\{|p|,{Q}\}^{2m}\dd x\right)^{1/2}
\left(\io |\nabla p|^2 \min\{|p|,{Q}\}^{2m}\dd x\right)^{1/2}\\
&\le& \frac{\mu_0}{2}(2{m}+1)\io |\nabla p|^2 \min\{|p|,{Q}\}^{2m}\dd x + \frac{2{m}+1}{2\mu_0} U^2(t)
\io \min\{|p|,{Q}\}^{2m}\dd x\,.
\eers
For the boundary term, we have
$$
\ipo \gamma(x) (p^* - p)h_{{Q},m}(p)\dd s(x)
\le \frac{1}{2{m}+2}\ipo\gamma(x)(|p^*|^{2{m}+2} -\min\{|p|,{Q}\}^{2{m}+2})\dd s(x)\,.
$$
On the left hand side of \eqref{m2}, we have
$$
\int_0^t f'(p)h_{{Q},m}(p) p_t(x,\tau)\dd \tau = F_{{Q},m}(p(x,t))- F_{{Q},m}(p(x,0))\,,
$$
where we set $F_{{Q},m}(p) = \int_0^p f'(z)h_{{Q},m}(z) \dd z$. We claim that for every $p\in \real$ we have
\be{pl2}
\frac{f_3}{2{m}+2} |p|^{2{m}} \ge F_{{Q},m}(p) \ge \frac{f_2}{4{m}} (\min\{|p|,{Q}\}^{2{m}} - 1)\,.
\ee
The upper bound is easy. We have for $z > 0$ that $f'(z)h_{{Q},m}(z) \le f_3 z |z|^{2m}$ and similarly
for $z<0$, and it suffices to integrate. To get the lower bound, set
$$
\hat F_{{Q},m}(p) = F_{{Q},m}(p) - \frac{f_2\min\{|p|,{Q}\}^{2{m}}}{4m}\,.
$$
Then for $p > {Q}$ we have $\hat F_{{Q},m}'(p)  = f'(p)h_{{Q},m}(p) > 0$. For $p \in (0,{Q})$ we have
$$
\hat F_{{Q},m}'(p)  = f'(p)p|p|^{2m} - \frac{f_2 p|p|^{2m-1}}{2} \ge \frac{f_2 p|p|^{2m-1}}{2(1+p^2)}(p^2 - 1),
$$
hence the minimum of $\hat F_{{Q},m}(p)$ is attained at $p=1$, with $\hat F_{{Q},m}(1) \ge - \frac{f_2}{4m}$,
which is exactly \eqref{pl2}. The case $p<0$ is symmetric.

Summarizing the above estimates, we obtain by integrating \eqref{m2} from $0$ to $t$ that
\ber\nonumber
&&\hspace{-10mm}\io\min\{|p|,{Q}\}^{2{m}}(x,t)\dd x
+ \frac{\mu_0}{f_2}2m(2m+1)\int_0^t \io |\nabla p|^2 \min\{|p|,{Q}\}^{2m}\dd x\dd\tau\\ \nonumber
&&+ \frac{2m}{f_2({m}+1)}\ipo\gamma(x)\min\{|p|,{Q}\}^{2{m}+2}\dd s(x)\\ \label{m4}
&\le& C(1+m^2) \int_0^t U^2(\tau)\left(K^{2m} + \io\min\{|p|,{Q}\}^{2{m}}(x,\tau)\dd x\right)\dd\tau\,.
\eer
In particular, the function $w_m(t) := \io\min\{|p|,{Q}\}^{2{m}}(x,t)\dd x$ satisfies the inequality
$$
w_m(t) \le C(1+m^2) \int_0^t U^2(\tau)\left(K^{2m} + w_m(\tau)\right)\dd\tau\,,
$$
and by Gronwall's argument (note that $U^2 \in L^1(0,T)$ by \eqref{m3}), there exists a constant $C(m)$
depending on $m$ and independent of ${Q}$ such that $\sup_{t\in [0,T]} w_m(t) \le C(m)$. Hence,
we can let ${Q}$ tend to $\infty$ in \eqref{m4} and obtain
\ber\nonumber
&&\hspace{-10mm}\io|p|^{2{m}}(x,t)\dd x
+ \frac{\mu_0}{f_2}2m(2m+1)\int_0^t \io |\nabla p|^2 |p|^{2m}\dd x\dd\tau\\ \nonumber
&&+ \frac{2m}{f_2({m}+1)}\ipo\gamma(x)|p|^{2{m}+2}\dd s(x)\\ \label{m4a}
&\le& C(1+m^2) \int_0^t U^2(\tau)\left(K^{2m} + \io|p|^{2{m}}(x,\tau)\dd x\right)\dd\tau\,.
\eer
In particular,
\be{m4b}
p \in L^\infty(0,T; L^q(\Omega))\quad \forall q \ge 1\,,
\ee
but the norm of $p$ in this space still depends on $q$. Note that
\ber\nonumber
\int_0^t \io |\nabla p|^2 |p|^{2m}\dd x\dd\tau &\ge& \int_0^t \io |\nabla p|^2 |p|^{2{m}-2}\dd x\dd\tau
- \int_0^t \io |\nabla p|^2\dd x\dd\tau \\ \label{cm3}
&\ge& \int_0^t \io |\nabla p|^2 |p|^{2{m}-2}\dd x\dd\tau - C
\eer
by virtue of \eqref{ee8}. Indeed, it suffices to split the integration domain into the parts where
$p\ge 1$ and $p<1$. Using \eqref{cm3}, we rewrite \eqref{m4} as
\ber\nonumber
&&\hspace{-10mm}\io|p(x,t)|^{2{m}}\dd x + \frac{\mu_0}{f_2}\frac{2(2m+1)}{m}
\int_0^t \io |\nabla (p|p|^{m-1})|^2 \dd x\dd\tau +
\frac{2m}{f_2({m}+1)}\ipo\gamma(x)|p|^{2{m}+2}\dd s(x)\\ \label{m5}
&\le& C(1+m^2) \int_0^t U^2(\tau)\left(K^{2m} + \io|p(x,\tau)|^{2{m}}\dd x\right)\dd\tau\,.
\eer
Putting $c:= 4\mu_0/f_2$, we obtain \eqref{m6}, and Lemma \ref{l1} is proved.
\epf

\bpf{Proof of Proposition \ref{p2}}
Set $w := p |p|^{m-1}$ for $p = p\oi$. Then \eqref{m6} reads
\be{m7}
|w(t)|_2^2 + c \int_0^t|\nabla w(\tau)|_2^2\dd\tau
\le C(1+m^2) \int_0^t U^2(\tau)\left(K^{2m} + |w(\tau)|_2^2\right)\dd\tau\,.
\ee
We now choose $s = 3/4$, and invoke the Gagliardo-Nirenberg inequality \eqref{gn} in the form
$$
|w(\tau)|_2 \le C\left(|w(\tau)|_{2s} + |w(\tau)|_{2s}^{1-\gamma} |\nabla w(\tau)|_2^\gamma\right)\,,
\quad \gamma = \frac{\frac{1}{2s} - \frac{1}{2}}{\frac{1}{2s} - \frac{1}{6}} = \frac13\,.
$$
H\"older's inequality enables us to rewrite \eqref{m7} as
\bers
&&\hspace{-10mm}|w(t)|_2^2 + c \int_0^t|\nabla w(\tau)|_2^2\dd\tau\\
&\le& C(1+m^2) \int_0^t U^2(\tau)\left(K^{2m} + |w(\tau)|_{2s}^2\right)\dd\tau
+ C(1+m^2)\int_0^t U^2(\tau)|w(\tau)|_{2s}^{4/3}|\nabla w(\tau)|_2^{2/3}\dd\tau
\\
&\le& C(1+m^2) \int_0^t U^2(\tau)\left(K^{2m} + |w(\tau)|_{2s}^2\right)\dd\tau\\
&& +\ C(1+m^2)\left(\int_0^t U^3(\tau)|w(\tau)|_{2s}^2\dd\tau\right)^{2/3}
\left(\int_0^t|\nabla w(\tau)|_2^2\dd\tau\right)^{1/3},
\eers
and we conclude that
\be{m8}
|w(t)|_2^2 \le C(1+m^3) \int_0^t U^3(\tau)\left(K^{2m} + |w(\tau)|_{2s}^2\right)\dd\tau\,,
\ee
or, in terms of $p$,
\be{m9}
|p(t)|_{2m}^{2m} \le C(1+m^3) \int_0^t U^3(\tau)\left(K^{2m} + |p(\tau)|_{2s m}^{2m}\right)\dd\tau\,.
\ee
We have $U \in L^3(0,T)$ by \eqref{m3}. If now $V$ is a constant such that
$\max\{K, |p(\tau)|_{2s m}\} \le V$, then
\be{m10}
\tilde V:= \max\{K, |p(t)|_{2m}\} \le \big(C(1+m^3)\big)^{1/2m} V\,.
\ee
We now define $m_k := s^{-k} = (4/3)^k$, and put
$$
V_k:= \max\{K, \sup_{\tau\in (0,T)}|p(\tau)|_{2m_k}\}\,.
$$
By \eqref{m10} we have for all $k=2, 3, \dots$ that
\be{m11}
V_k \le \left(C\left(1 + \left(\frac43\right)^{3k}\right)\right)^{((3/4)^k)/2} V_{k-1}\,,
\ee
hence
\be{m12}
\log V_k - \log V_{k-1} \le \frac12\left(\frac34\right)^{k}
\log\left(C\left(1 + \left(\frac43\right)^{3k}\right)\right).
\ee
The right hand side of \eqref{m12} is a convergent series, and $V_1$ is finite by virtue
of \eqref{m4b}, so that we can conclude that the sequence $\{V_k\}$ is uniformly bounded
independently of $i$, which we wanted to prove.
\epf


\subsection{Estimate 11}\label{est6}

Test \eqref{e24r} by $\psi = M(p)_t$, with $M(p)$ given by \eqref{mum}.
From \eqref{mono} and Proposition \ref{p2}, it follows that there exists
a constant $c>0$ such that for every $t \in (0,T)$ we have
\ber \nonumber 
&&\hspace{-9mm} c\io p_t^2(x,t) \dd x + \frac{1}{2\rho_L}\frac{\dd}{\dd t}\io |\nabla M(p)|^2(x,t) \dd x +
\frac{\dd}{\dd t}\ipo \gamma_p(x) \hat M(p)(x,t)\dd s(x)\\ \nonumber
&&\le \io (|\dive u_t| |M(p)_t|)(x,t)\dd x
+ \frac{\dd}{\dd t} \ipo \gamma_p(x) (M(p) p^*)(x,t)\dd s(x)\\ \label{e61}
&&\quad - \ipo\gamma_p(x) (M(p) p^*_t)(x,t)\dd s(x)
\eer
with $\hat M(p) = \int_0^p p'\mu(p')\dd p'$. Hence, using also \eqref{ee10}, we have for all $t\in (0,T)$ that
\be{e62}
\|p\oi_t\|_2^2 + |\nabla p\oi(t)|_2^2 + \ipo \gamma_p(x) |p\oi|^2(x,t)\dd s(x) \le C\,.
\ee
By comparison in Eq.~\eqref{e24r}, we also obtain that
\be{e63}
\|M(p\oi)\|_{L^2(0,T; W^{2,2}(\Omega))} \le C\,.
\ee


\subsection{Estimate 12}\label{est7}

By the same argument based on \eqref{e62}--\eqref{e63}, we have similarly as in \eqref{ee14} that
\be{e701}
\|\nabla p\oi\|_{L^p(0,T;L^q(\Omega))} \le C \ \for q = \frac{3}{1 - 3\alpha},\ p = \frac{4}{6\alpha + 1}
\ee
for every $\alpha \in [0, 1/6)$. In particular,
for $\alpha = 1/30$, we have $\partial_{x_i} u_t, \partial_{x_i} p
\in L^{10/3}(\Omega\times (0,T))$. Hence, Eq.~\eqref{e25r} has the form
\be{e25t}
\io \big(c_0 \theta_t \zeta
 + \kappa(\theta) \nabla\theta\cdot\nabla\zeta\big) \dd x 
 + \ipo \gamma_\theta(x)(\theta - \theta^*)\zeta\dd s(x)
= \io(A(x,t) + B(x,t) \theta)\zeta\dd x
\ee
with $A \in L^{5/3}(\Omega\times (0,T))$, $B \in L^{10/3}(\Omega\times (0,T))$ bounded independently of $i$.
We test \eqref{e25t} with $\zeta = \theta^r$ (with $r$ to be chosen later) and obtain using H\"older's inequality for every $t\in (0,T)$ that
\ber \nonumber \label{e71a}
&&\hspace{-16mm}\frac{1}{r+1}\io \theta^{r+1}(x,t)\dd x
+ r\int_0^T\!\io \theta^{r-1}\kappa(\theta)|\nabla\theta|^2(x,\tau)\dd x\dd\tau
+ \int_0^T\!\ipo \gamma_\theta(x)\theta^{r+1}(x,\tau)\dd s(x)\dd \tau\\
&\le& C\left(1 + \|\theta\|_{5r/2}^r + \|\theta\|_{10(r+1)/7}^{r+1}\right)
\eer
with, by Hypothesis \ref{h1}\,((ii),
\be{e71}
\theta^{r-1}\kappa(\theta)|\nabla\theta|^2\ge \frac{1}{C}\theta^{r+a}|\nabla\theta|^2.
\ee
We already have the estimate \eqref{ee6}. Assume that for some $z \ge 8/3$ we have proved
\be{e72}
\|\theta\|_z \le C_0(z)
\ee
with some $C_0(z)>0$. For this value of $z$ we choose in \eqref{e71a}
\be{e73}
r = r(z) = \left\{
\begin{array}{ll}
(7z/10) - 1 & \for z \in [8/3, 10/3]\,,\\
2z/5 & \for z > 10/3\,.
\end{array}
\right.
\ee
Then $\|\theta\|_{5r/2}^r + \|\theta\|_{10(r+1)/7}^{r+1} \le C(1+ \|\theta\|_z^{r+1})$,
and we have by virtue of \eqref{e71a}--\eqref{e71} that
\be{e74}
\frac{1}{r+1}\io \theta^{r+1}(x,t)\dd x + r\int_0^T\io \theta^{r+a}|\nabla\theta|^2(x,\tau)\dd x\dd\tau
\le C\left(1 + C_0(z)^{r+1}\right).
\ee
Set
\be{e75}
p = \frac{r+a}{2} + 1\,, \ s = \frac{r+1}{p}\,, \ w = \theta^p\,.
\ee
Then \eqref{e74} can be written as
\be{e76}
\frac{1}{r+1}\io w^s(x,t)\dd x + r\int_0^T\io |\nabla w|^2(x,\tau)\dd x\dd\tau \le C\left(1 + C_0(z)^{r+1}\right).
\ee
For $s<q<6$ we have by virtue of the Gagliardo-Nirenberg inequality \eqref{gn} that
\be{e77}
|w(\tau)|_q \le C\left(|w(\tau)|_{s} + |w(\tau)|_{s}^{1-\gamma} |\nabla w(\tau)|_2^\gamma\right)\,,
\quad \gamma = \frac{\frac{1}{s} - \frac{1}{q}}{\frac{1}{s} - \frac{1}{6}}\,.
\ee
If $q$ is chosen in such a way that $q\gamma = 2$, that is,
\be{e78}
q = \frac23 s + 2\,,
\ee
then it follows from \eqref{e77} and Young's inequality that
\be{e77e}
\|w\|_q \le C\left(\sup_{\tau\in[0,T]}|w(\tau)|_{s} + \sup_{\tau\in[0,T]}|w(\tau)|_{s}^{1-\gamma}
\|\nabla w\|_2^{2/q}\right)
\le C\left(\sup_{\tau\in[0,T]}|w(\tau)|_{s} + \|\nabla w\|_2\right)\,.
\ee
By \eqref{e76}, we have
\bers
\sup_{\tau\in[0,T]}|w(\tau)|_{s} &\le& C\left((1+r)\big(1 + C_0(z)^{r+1}\big)\right)^{1/s}\,,\\
\|\nabla w\|_2 &\le& C\big(1 + C_0(z)^{r+1}\big)^{1/2} \le C\big(1 + C_0(z)^{r+1}\big)^{1/s}
\eers
(note that $s<2$), so that
\be{e79}
\|w\|_q \le C\left((1+r)\big(1 + C_0(z)^{r+1}\big)\right)^{1/s}\,.
\ee
Putting
$$
\hat z = pq = \frac53 r(z) + \frac83 + a\,,
$$
we have by \eqref{e79} that
\be{e79a}
\|\theta\|_{\hat z} = \|w\|_q^{1/p} \le C\left((1+r)\big(1 + C_0(z)^{r+1}\big)\right)^{1/(ps)}\,,
\ee
and using the identity $ps = r+1$, we obtain the implication
\be{e710}
\|\theta\|_z \le C_0(z) \Longrightarrow \|\theta\|_{\hat z} \le \hat C_0(z)\,,\quad
\hat C_0(z) = C(1+ C_0(z))\,.
\ee
The sequence
$$
z_k = \frac53 r(z_{k-1}) + \frac83 + a\,, \quad z_0 = \frac83
$$
converges to $z_\infty = 8 + 3a$. After finitely many iterations we obtain
\be{e711}
\|\theta\|_z \le C \quad \mbox{for every }\ z < 8 + 3a\,.
\ee
Consequently, by \eqref{e71a}--\eqref{e74}, we have for $r < \frac{16}{5} + \frac65$ and $t\in (0,T)$ that
\be{e712}
\io \theta^{r+1}(x,t)\dd x + r\int_0^T\!\io \theta^{r-1}\kappa(\theta)|\nabla\theta|^2(x,\tau)\dd x\dd\tau
+ \int_0^T\!\ipo \gamma_\theta \theta^{r+1}(x,\tau)\dd s(x)\dd \tau \le C\,. 
\ee
This enables us to derive an upper bound for the integral
$\io \kappa(\theta) \nabla\theta\cdot\nabla\zeta \dd x$
which we need for getting an estimate for $\theta_t$ from the equation \eqref{e25t}. We have
by H\"older's inequality and Hypothesis \ref{h1}\,(ii) that
\ber \nonumber
\io |\kappa(\theta) \nabla\theta\cdot\nabla\zeta| \dd x &=&
\io |\theta^{(1-r)/2}\kappa^{1/2}(\theta) \theta^{(r-1)/2}\kappa^{1/2}(\theta)\nabla\theta\cdot\nabla\zeta| \dd x\\
\label{e715}
&\le& C\left(\io \theta^{r-1}\kappa(\theta) |\nabla\theta|^2\dd x\right)^{1/2}
 \left(\io \theta^{2+b-r}|\nabla\zeta|^2\dd x\right)^{1/2}.
\eer
We now want to choose $\eta > 0$ and $r$ satisfying \eqref{e712} such that
\be{e714}
2+b-r = \frac{r+1}{1+\eta}\,,
\ee
or, equivalently,
\be{e713}
r = \frac{(1+2\eta) + b (1+\eta)}{2+\eta} = \frac{1+b}{2} + \frac{\eta(b+3)}{2(2+\eta)}\,.
\ee
If $b$ satisfies the condition in Hypothesis \ref{h1}\,(ii), then there exists $\eta \in (0,1)$
such that for $r$ given by \eqref{e713} the condition in \eqref{e712} holds. Hence, by \eqref{e715},
\be{e716}
\io |\kappa(\theta) \nabla\theta\cdot\nabla\zeta| \dd x 
\le C\left(\io \theta^{r-1}\kappa(\theta) |\nabla\theta|^2\dd x\right)^{1/2}
 \left(\io \theta^{r+1} \dd x\right)^{1/(1+\eta)}\left(\io |\nabla\zeta|^{q^*}\dd x\right)^{1/q^*}
\ee
with $q^* = \frac{2(1+\eta)}{\eta}$. Hence, by \eqref{e712},
\be{e717}
\int_0^T\io |\kappa(\theta\oi) \nabla\theta\oi\cdot\nabla\zeta| \dd x \dd t
\le C \|\zeta\|_{L^2(0,T;W^{1,q^*}(\Omega))}\,.
\ee
By \eqref{ee14} and \eqref{e701} for $\alpha = 0$, the functions $A$ and $B$ in \eqref{e25t}
satisfy uniform bounds $A \in L^2(0,T; L^{3/2}(\Omega))$, $B \in L^4(0,T; L^{3}(\Omega))$
independent of $i$, so that testing with $\zeta \in L^2(0,T;W^{1,q^*}(\Omega))$ is admissible.
We thus obtain from \eqref{e25t} that
\be{e718}
\int_0^T\io \theta\oi_t \zeta \dd x \dd t \le C \|\zeta\|_{L^2(0,T;W^{1,q^*}(\Omega))}\,.
\ee


\subsection{Passage to the limit as $i \to \infty$}\label{iinf}

In the system \eqref{e23r}--\eqref{e25r} with $\delta = \delta_i$, $R = R_i$,
and $(u, p, \theta) = (u\oi, p\oi, \theta\oi)$, we fix test functions
$\phi \in W^{2,2}(\Omega; \real ^3) \cap X_0$, $\psi \in X$, and $\zeta \in L^2(0,T;W^{1,q^*}(\Omega))$
with $q^*$ from \eqref{e717}.
The term $\delta_i \Bcal u\oi_{tt}$ in \eqref{e23r} converges to $0$ in $L^2$ by \eqref{ee10},
the regularization $K_{R_i}(p)$ vanishes for $R_i > C^*$ by Proposition \ref{p2}.

By \eqref{ee10} and \eqref{ee11}, the sequence $\{\nas u\oi_t\}$ is precompact in $L^4(0,T; L^q(\Omega; \tens))$
for every $q \in [1,3)$. Similarly, by \eqref{e62}--\eqref{e63}, $\{\nabla p\oi\}$ is precompact in
$L^4(0,T; L^q(\Omega; \real^3))$ for every $q \in [1,3)$, and $\{p\oi\}$ is precompact in
$L^q(\Omega; C[0,T])$ for every $q > 1$ by virtue of Proposition \ref{p2}. Finally, by \eqref{ee5},
\eqref{e711}, and \eqref{e718}, $\{\theta\oi\}$ is precompact e.\,g. in $L^8(\Omega\times (0,T))$.
Hence, using also \eqref{e717}, we select a subsequence and pass to the weak limit
in the linear terms in \eqref{e23r}--\eqref{e25r}, and to the strong limit in all
nonlinear non-hysteretic terms. Obviously, if $\theta\oi$
converge strongly to $\theta$ in $L^8(\Omega\times (0,T))$, then $Q_{R_i}(\theta\oi) \to \theta$ strongly
in $L^8(\Omega\times (0,T))$ as well, and if $|\nabla p\oi|^2 \to |\nabla p|^2$
in $L^2(0,T; L^{q/2}(\Omega))$ strongly, then 
$Q_{R_i}(|\nabla p\oi|^2) \to |\nabla p|^2$ in $L^2(0,T; L^{q/2}(\Omega))$ strongly.
By the same argument as in Subsection \ref{ninf}, we show that the hysteresis terms
$G[p\oi]_t$, $\|D_P[\nas u\oi]_t\|_*$, $|D_G[p\oi]_t|$ converge weakly in $L^2(\Omega\times (0,T))$,
and that the limit as $i \to \infty$ yields a solution to \eqref{e23v}--\eqref{e25v}
with $\phi \in W^{2,2}(\Omega; \real^3) \cap X_0$. By density we conclude that $\phi \in X_0$
is an admissible test function, which completes the proof of Theorem~\ref{t1}.

\end{document}